\documentclass{article}

\usepackage[english]{babel}
\usepackage[UTF8]{ctex}
\usepackage[letterpaper,top=2cm,bottom=2cm,left=3cm,right=3cm,marginparwidth=1.75cm]{geometry}
\usepackage{tikz}
\usetikzlibrary{intersections}
\usetikzlibrary{arrows.meta,decorations,calligraphy}
\usepackage{amsmath}
\usepackage{graphicx}
\usepackage[colorlinks=true, allcolors=blue]{hyperref}
\usepackage{graphicx}
\usepackage{amsfonts}
\usepackage{amstext}

\usepackage{amsthm}
\usepackage{amssymb,amsmath,amsbsy}
\usepackage{verbatim}
\usepackage{thm-restate}
\newcommand{\p}{\mathcal{P}}
\newcommand{\dff}{\draw[fill,color=gray!10]}
\newcommand{\df}{\draw[fill]}

\newcommand{\bt}{\begin{tikzpicture}}
\newcommand{\et}{\end{tikzpicture}}
\providecommand{\keywords}[1]
{
  \small	
  \textbf{\textit{Keywords---}} #1
}
\newtheorem{theorem}{Theorem}

\newtheorem{corollary}[theorem]{Corollary}
\newtheorem{lemma}[theorem]{Lemma}
\newtheorem{definition}[theorem]{Definition}

\tikzset{dot/.style={font=\LARGE}}
\title{The planar Tur\'an number of $\{K_4,C_5\}$ and $\{K_4,C_6\}$}

\author{Ervin Gy\H{o}ri\thanks
  {R\'enyi Institute, Budapest,
 Hungary.  Research partially supported by the NKFIH Grant 132696. E-mail: {\tt gyori.ervin@renyi.hu}}
  \\ Alan Li\thanks{ Amherst College, 220 South Pleasant Street,
Amherst, MA 01002  E-mail:
 {\tt ali24@amherst.edu}}
 \qquad
Runtian Zhou \thanks{
Davidson College, Davidson, North Carolina 28035
Email: {\tt dazhou@davidson.edu}}
}

\begin{document}
\maketitle

\begin{abstract}
Let $\mathcal{H}$ be a set of graphs. The planar Tur\'an number, $ex_\p(n,\mathcal{H})$, is the maximum number of edges in an $n$-vertex planar graph which does not contain any member of $\mathcal{H}$ as a subgraph. When $\mathcal{H}=\{H\}$ has only one element, we usually write $ex_\p(n,H)$ instead. The topic of extremal planar graphs was initiated by Dowden (2016). He obtained sharp upper bound for both $ex_\p(n,C_5)$ and $ex_\p(n,K_4)$. Later on, we obtained sharper bound for $ex_\p(n,\{K_4,C_7\})$. In this paper, we give upper bounds of $ex_\p(n,\{K_4,C_5\})\leq {15\over 7}(n-2)$ and $ex_\p(n,\{K_4,C_6\})\leq {7\over 3}(n-2)$. We also give constructions which show the bounds are tight for infinitely many graphs.
\end{abstract}
\section{Introduction and Main Results}

In this paper, all graphs considered are planar, undirected, finite and contain neither loops nor multiple edges (unless otherwise stated). Actually, we study plane graphs what are embeddings (drawings) of graphs in the plane such that the edge curves do not cross each other, they may share just endpoints of the curves. We use $C_k$ to denote the cycle of $k$ vertices and $K_r$ to denote the complete graph of $r$ vertices. We use $n$-face to denote a face with $n$ edges.

The Tur\'an number $ex(n,H)$ for a graph $H$ is the maximum number of edges in an $n$-vertex graph with no copy of $H$ as a subgraph. The first result on the topic of Tur\'an number was obtained by Mantel and Tur\'an, who proved that the balanced complete $r$-partite graph is the unique extremal graph of $ex(n,K_{r+1})$ edges. The Erd\H{o}s-Stone-Simonovits \cite{erdos1963structure,erdos1962number} then generalized this result and asymptotically determined $ex(n,H)$ for all nonbipartite graphs $H:$ $ex(n,H)=(1-{1\over \mathcal{X}(H)-1}){n\choose 2}+o(n^2)$.

In 2016, Dowden et al. \cite{dowden2016extremal} initiated the study of Tur\'an-type problems when host graphs are plane, i.e., how many edges can a plane graph on $n$ vertices have, without containing a given smaller graph? Let $\mathcal{H}$ be a set of graphs. The planar Tur\'an number, $ex_\p(n,\mathcal{H})$, is the maximum number of edges in an $n$-vertex planar graph which does not contain any member of $\mathcal{H}$ as a subgraph. When $\mathcal{H}=\{H\}$ has only one element, we usually write $ex_\p(n,H)$ instead. Dowden et al. \cite{dowden2016extremal} obtained the tight bounds $ex_{\mathcal{P}}(n,K_4) =3n-6$, for all $n\geq 4$, $ex_{\mathcal{P}}(n,C_4) \leq\frac{15(n-2)}{7}$, for all $n\geq 4$ and $ex_{\mathcal{P}}(n,C_5) \leq\frac{12n-33}{5}$, for all $n\geq 11$. For $k\in \{4,5\}$, let $\Theta_k$ denote the graph obtained from $C_k$ by adding a chord. Y. Lan et al. \cite{lan2019extremal} showed that $ex_\p(n,\Theta_4)\leq {15(n-2)\over 5}$ for all $n\geq 4$, $ex_\p(n,\Theta_5)\leq {5(n-2)\over 2}$ for all $n\geq 5$ and $ex_\p(n,C_6)\leq {18(n-2)\over 7}$. The bounds for $ex_\p(n,\Theta_4)$ and $ex_\p(n,\Theta_5)$ are tight infinitely often, but the upper bound for $ex_\p(n,C_6)$ was improved by D. Ghosh et al. \cite{ghosh2022planar}. They proved $ex_\p(n,C_6)\leq {5n-14\over 2}$ for all $n\geq 18$. Recently, R. Shi et al. \cite{shi2023planar} and us \cite{gyori2023planar} independently proved sharp bound for $ex_\p(n,C_7)\leq {18\over 7}n-{48\over 7}$ for all $n\geq 60$. In this paper, we are interested in $ex_\p(n,\{K_4,C_m\})$ where $m$ is an integer. The first result of such problems is $ex_\p(n,\{K_4,C_4\})=ex_{\mathcal{P}}(n,K_4) =3n-6$ as was proved by Dowden et al. \cite{dowden2016extremal}. In \cite{gyori2023planar}, we proved $ex_\p(n,\{K_4,C_7\})\leq {18\over 7}n-{48\over 7}$ for all $n\geq 60$ and showed the bound is tight for infinitely many graphs. In this paper, we focus on $ex_\p(n,\{K_4,C_5\})$ and $ex_\p(n,\{K_4,C_6)$ and give bounds that are sharp for infinitely many graphs.

We denote the vertex and the edge sets of a graph $G$ by $V(G)$ and $E(G)$ respectively. We also denote the number of vertices and edges of $G$ by $v(G)$ and $e(G)$ respectively. The minimum degree of $G$ is denoted $\delta(G)$. The main result is based on the theorem as follows:

\begin{theorem}\label{t1}
Let $G$ be a $\{K_4,C_5\}$-free plane graph on $n$ $(n\geq 5)$ vertices with $\delta(G)\geq 3$. Then $$e(G)\leq {15\over 7}(n-2).$$
\end{theorem}
We use Theorem \ref{t1} in order to establish our desired result, which gives the upper bound of ${15\over 7}(n-2)$ for all $\{K_4,C_5\}$-free plane graphs with at least $30$ vertices.

\begin{corollary}\label{c1}
Let $G$ be a $\{K_4,C_5\}$-free plane graph on $n$ $(n\geq 15)$ vertices. Then $$e(G)\leq {15\over 7}(n-2).$$
\end{corollary}

\begin{corollary}\label{c12}
Let $G$ be a $\{K_4,C_5\}$-free plane graph on $n$ $(n\geq 2)$ vertices such that $G$ contains no maximal $2$-connected subgraph of more than $4$ vertices. Then $$e(G)\leq 2n-{15\over 7}.$$
\end{corollary}

We show that Corollary \ref{c1} is tight for infinitely many  graphs.

\begin{theorem}\label{t2}
For every $n\geq 9$, $n\equiv 2\mod 7$, there exists a $\{K_4,C_5\}$-free plane graph $G$ with $e(G)= {15\over 7}(n-2).$
\end{theorem}

\begin{theorem}\label{t3}
Let $G$ be a $\{K_4,C_6\}$-free plane graph on $n$ $(n\geq 6)$ vertices with $\delta(G)\geq 3$. Then $$e(G)\leq {7\over 3}(n-2).$$
\end{theorem}
We use Theorem \ref{t3} in order to establish our desired result, which gives the upper bound of ${7\over 3}(n-2)$ for all $\{K_4,C_6\}$-free plane graphs with at least $9$ vertices.

\begin{corollary}\label{c2}
Let $G$ be a $\{K_4,C_6\}$-free plane graph on $n$ $(n\geq 9)$ vertices. Then $$e(G)\leq {7\over 3}(n-2).$$
\end{corollary}
\begin{corollary}\label{c22}
Let $G$ be a $\{K_4,C_6\}$-free plane graph on $n$ $(n\geq 2)$ vertices such that $G$ contains no maximal $2$-connected subgraph of more than $5$ vertices. Then $$e(G)\leq {31\over 15}n-{7\over 3}.$$
\end{corollary}
We show that  Corollary \ref{c2} is tight for infinitely many graphs.

\begin{theorem}\label{t4}
For every $n\equiv 50\mod 288$, there exists a $\{K_4,C_6\}$-free planar graph $G$ with $e(G)= {7\over 3}(n-2).$
\end{theorem}

\section{Definitions}
\begin{definition}\label{tb}
Let $G$ be a plane graph. Then a \textbf{triangular block} $B$ is a subset of the edge set of $G$. Triangular blocks are built as follows: \\
1) Begin with an edge $e \in E(G)$. If it is not in any $3$-face of $G$, then we have $B= \{e\}$; \\
2) Otherwise, we add $e$ into $B$, and for each edge $3$-face containing some edge of $B$, we add all edges of that face into $B$; \\
3) Repeat step 2) until no more edges can be added into $B$.

Similar to \cite{ghosh2022planar}, we make the following observations: \\
(i) Given a triangular block $B$, no matter which edge we begin with, we always obtain $B$ as the triangular block. \\
(ii) The triangular blocks form a partition of the edges in $G$.

Thus, given a plane graph $G$, the set of its triangular blocks is well defined. Denote this set as $B(G)$.
\end{definition}

\begin{definition}
Let $G$ be a plane graph and $e\in E(G)$. Let the two faces incident to $e$ have length $l_1$ and $l_2$. If $e$ is incident to only one face, let $l_1=l_2$ be the length of that face. The contribution of $e$ to the face number of $G$, denoted by $f(e)$, is defined as $$f(e)={1\over l_1}+{1\over l_2}.$$
Let $B\subseteq E(G)$ be a triangular block. The contribution of $B$ to the face number of $G$, denoted as $f(B)$, is defined as $$f(B)=\sum_{e\in B}f(e).$$The contribution of $B$ to the edge number of $G$, denoted as $e(B)$, is defined as the size of $B$.
\end{definition}

\section{Proof of Theorem \ref{t2}: Extremal Graph Construction}
\begin{figure}
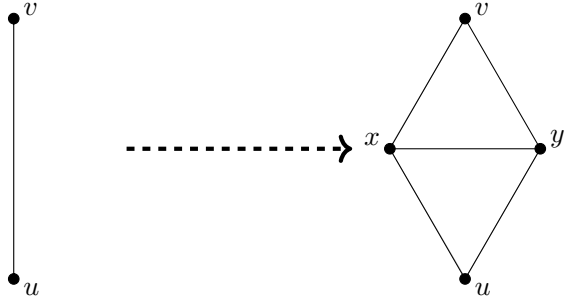

\centering
\bt
\df (-3,1.732)circle(2pt)--(-3,-1.732)circle(2pt);
\df (3,1.732)circle(2pt)--(2,0)circle(2pt)--(3,-1.732)circle(2pt)--(4,0)circle(2pt);
\draw (3,1.732)--(4,0)--(2,0);
\draw[dashed,ultra thick,->] (-1.5,0)--(1.5,0);
\draw (-3,1.732)node[anchor=-150]{$v$};
\draw (-3,-1.732)node[anchor=150]{$u$};
\draw (3,1.732)node[anchor=-150]{$v$};
\draw (3,-1.732)node[anchor=150]{$u$};
\draw (2,0)node[anchor=-30]{$x$};
\draw (4,0)node[anchor=-150]{$y$};
\et
\caption{\label{c5t}Local Transformation}
\end{figure}

\begin{proof}
This construction comes from the intuition that $B_4$(Figure \ref{b4}) is the ``optimal" triangular block, which will be implicitly implied throughout the proof of Theorem \ref{t1} as we achieve $15f(B_4)-8e(B_4)=0$.

Notice that for each $n'\geq 3$, we can find a planar triangulated graph $G'$ with $n'$ vertices. We than do local change to each edge of $G'$ as described in Figure \ref{c5t} to get graph $G$. Specifically, for each edge $(vu)\in E(G)$, we delete $(vu)$ and add two new vertices $x,y$ together with $5$ edges $(vx),(vy),(ux),(uy),(xy)$. Clearly $G$ is $\{K_4,C_5\}$-free. Denote $e(G')=e',v(G)=n,e(G)=e$. Because $G'$ is triangulated and planar, we know $e'=3n'-6$. Then,
\begin{align*}
e&=5e'=15n'-30\\
n&=n'+2e'=7n'-12\\
e&={15\over 7}(n-2).
\end{align*}
\end{proof}
\section{Proof of Theorem \ref{t1}}
\begin{proof}
Let $G$ be a $\{K_4,C_5\}$-free plane graph on $n$ $(n\geq 5)$ vertices with $\delta(G)\geq 3$. Denote the number of faces in $G$ as $f(G)$. By Euler's formula, $e(G)\leq {15\over 7}(n-2)$ is equivalent with $15f(G)-8e(G)\leq 0$. Observe that $15f(G)-8e(G)=\sum_{B\in B(G)}(15f(B)-8e(B))$. Hence it suffices to prove $15f(B)-8e(B)\leq 0$ for each $B\in B(G)$. We will do case analysis on different $B$'s. Since $G$ is $\{K_4,C_5\}$-free, there are only $3$ possible different types of $B$.

\textbf{Case 1: }$B=K_2$.

Since $\delta(G)\geq 3$, this edge must be incident to two faces each with length at least $4$. Thus $$15f(B)-8e(B)\leq 15\cdot ({1\over 4}+{1\over 4})-8=-0.5<0.$$

\textbf{Case 2: }$B=K_3$. 

Because $\delta(G)\geq 3$, no edge of this $K_3$ can be adjacent to a face of length $4$. Thus $$15f(B)-8e(B)\leq 15\cdot (1+3\cdot {1\over 6})-8\cdot 3=-1.5<0.$$

\textbf{Case 3: }$B=B_4$.

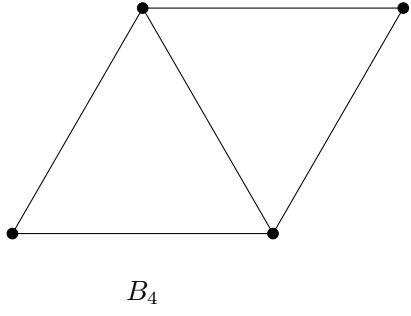
\begin{figure}
\centering
\begin{tikzpicture}
\draw[fill](210:2)circle(2pt);
\draw[fill](90:2)circle(2pt);
\draw[fill](30:4)circle(2pt);
\draw[fill](-30:2)circle(2pt);
\draw (210:2)--(90:2)--(30:4)--(-30:2)--cycle;
\draw(90:2)--(-30:2);
\draw (0,-1.5) node[below] {$B_{4}$};
\end{tikzpicture}
\caption{\label{b4}Triangular block on four vertices.}
\end{figure}
Since $G$ has no $C_5$ or $K_4$ and $\delta(G)\geq 3$, we know $$15f(B)-8e(B)\leq 15\cdot (2+4\cdot {1\over 6})-8\cdot 5=0.$$

\end{proof}

\section{Proof of Corollary \ref{c1}}
\begin{proof}
Let $G$ be a $\{K_4,C_5\}$-free plane graph on $n$ vertices. For simplicity of notations, in this section we denote $e(G)$ as $e$. We will show that either $e(G)\leq {15\over 7}(n-2)$, which is equivalent with $15n-7e\geq 30$, or $n<15$. We will repeat the following operation:

$1)$ Delete $x$ s.t. $deg(x)\leq 2$;

When we no longer can do this operation, we get an induced subgraph $G'$ with $\delta(G')\geq 3$. Denote 	$|E(G')|=e',|V(G')|=n'$.
Each time we do operation $1)$, denote the number of deleted edge as $e_d$. Observe that $15\cdot 1-7e_d\geq 1$. Thus we get $$15n-7e\geq 15n'-7e'+1\cdot(n-n').$$

In line with usual graph theoretic terminology, we call a maximal $2$-connected subgraph a \textbf{block}. Let $b$ be the total number of blocks of $G'$. Specifically, let $b_2,b_3$, and $b_4$ denote the number of blocks of size $2,3$, and $4$, respectively. Let $b_5$ denote the number of blocks of size at least $5$. Then we have $b=\sum_{i=2}^5b_i$.

Let's find a lower bound of $15n-7e-15$ for blocks of size $n$. Since $n$ is fixed, we need to plug in the largest possible $e$, which ideally would be the number of edges in a triangulated planar graph.

When $n\geq 5$, by Theorem \ref{t1} we have $15n-7e-15\geq 30-15=15$.

When $n=4$, since there is no $K_4$ as subgraph of $G$, we have $15n-7e-15\geq 15 \cdot 4-7\cdot5-15=10$ as in $B_4$.

When $n= 3$, as in $K_{3}$ we have $15n-7e-15\geq 15\cdot3-7\cdot3-15=9$.

When $n= 2$, as in $K_{2}$ we have $15n-7e-15\geq 15\cdot2-7\cdot1-15=8$.
\\

Denote the vertex number of the $i$th block as $n_i$ and edge number as $e_i$. Using the lower bounds presented above, when $n>n'$, we have
\begin{align*}
15n'-7e'&\geq15(\sum_{i=1}^{b}n_i-(b-1))-7\sum_{i=1}^{b}e_i\\
&=\sum_{i=1}^b(15n_i-7e_i-15)+15\\
&\geq 15b_5+10b_4+9b_3+8b_2+15\\
15n-7e&\geq 15b_5+10b_4+9b_3+8b_2+15+(n-n').
\end{align*}

If $b_5\geq 1$, then righthand side is at least $30$, as desired.

So, let us assume that $b_5=0$. It follows that
\begin{align*}
b&=\sum_{i=2}^4b_i\\
n&\leq \sum_{i=2}^4ib_i+(n-n')\\
15n-7e&\geq 15b_5+10b_4+9b_3+8b_2+15+(n-n')\\
&\geq 15+\sum_{i=2}^4ib_i+(n-n')\\
&\geq 15+n.
\end{align*}
This is greater than or equal to $30$ if $n\geq 15$.

Finally, if $n=n'$ and $n\geq 2$, since the last two vertices we deleted have at most one edge between them, we know $$15n-7e\geq 15\cdot 2-7\cdot 1+(n-2)=23+(n-2)$$which is greater than or equal to $30$ when $n\geq 9$. Since $9<15$, we complete the proof of Corollary \ref{c1}.

\end{proof}

\section{Proof of Theorem \ref{t4}: Extremal Graph Construction}
\begin{figure}
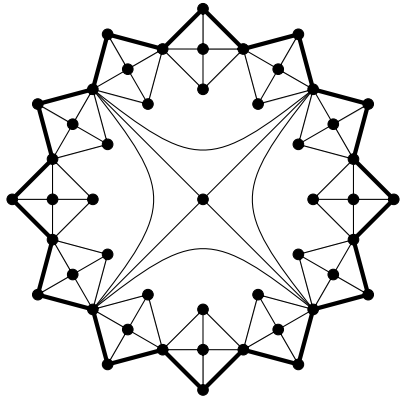

\centering
\bt
\foreach \x in {0,30,60,90,120,150,180,210,240,270,300,330}
{
\df (\x+15:2.071)--(\x:1.464)circle(2pt)--(\x-15:2.071)circle(2pt)--(\x:2)circle(2pt)--(\x:2.536)circle(2pt);
\draw[ultra thick](\x+15:2.071)--(\x:2.536)--(\x-15:2.071);
\draw (\x+15:2.071)--(\x:2)--(\x:1.464);
}
\foreach \x in {45,135,225,315}\draw (\x:2.071)..controls(\x+30:0.4)and(\x+60:0.4)..(\x+90:2.071);
\df (0,0)circle(2pt);
\draw (45:2.071)--(225:2.071);
\draw (135:2.071)--(315:2.071);
\et
\caption{\label{h0}The base case $H_0$}
\end{figure}
\begin{figure}
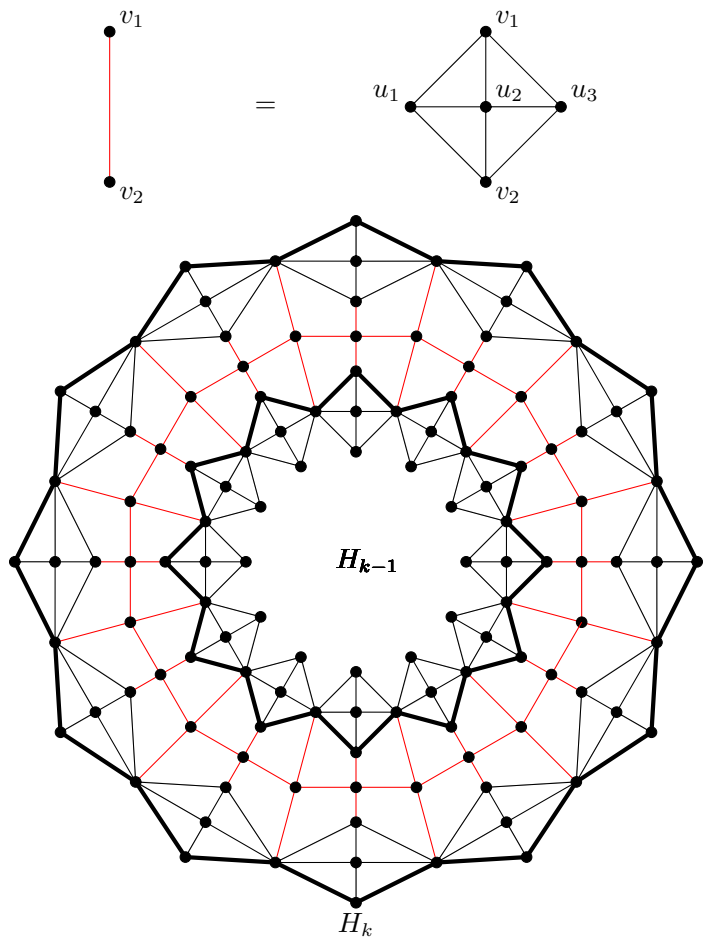

\centering
\bt
\draw[color=red](-2.5,-1)--(-2.5,1);
\df (-2.5,1)node[anchor=210]{$v_1$}circle(2pt);
\df (-2.5,-1)node[anchor=150]{$v_2$}circle(2pt);
\df (2.5,1)node[anchor=210]{$v_1$}circle(2pt)--(1.5,0)node[anchor=-30]{$u_1$}circle(2pt)--(2.5,-1)node[anchor=150]{$v_2$}circle(2pt)--(2.5,0)node[anchor=-150]{$u_2$}circle(2pt)--(3.5,0)node[anchor=210]{$u_3$}circle(2pt);
\draw(1.5,0)--(2.5,0)--(2.5,1)--(3.5,0)--(2.5,-1);
\draw(-0.7,0)node[anchor=180]{$=$};
\et

\bt

\foreach \x in {0,30,60,90,120,150,180,210,240,270,300,330}
{
\draw[color=red](\x-15:4.142)--(\x-15:3.106)--(\x-15:2.071);
\draw[color=red](\x:2.536)--(\x:3)--(\x:3.464);
\draw[color=red](\x+15:3.106)--(\x:3)--(\x-15:3.106);
\df(\x:3)circle(2pt);
\df(\x-15:3.106)circle(2pt);
\df (\x+15:4.142)--(\x:3.464)circle(2pt)--(\x-15:4.142)circle(2pt)--(\x:4)circle(2pt)--(\x:4.536)circle(2pt);
\draw[ultra thick](\x+15:4.142)--(\x:4.536)--(\x-15:4.142);
\draw (\x+15:4.142)--(\x:4)--(\x:3.464);
\df (\x+15:2.071)--(\x:1.464)circle(2pt)--(\x-15:2.071)circle(2pt)--(\x:2)circle(2pt)--(\x:2.536)circle(2pt);
\draw[ultra thick](\x+15:2.071)--(\x:2.536)--(\x-15:2.071);
\draw (\x+15:2.071)--(\x:2)--(\x:1.464);
\draw (-0.4,0)node[anchor=180]{$H_{k-1}$};
}
\draw (0,-4.536)node[anchor=90]{$H_k$};

\et
\caption{\label{hk}$H_k$ from $H_{k-1}$}
\end{figure}
\begin{figure}
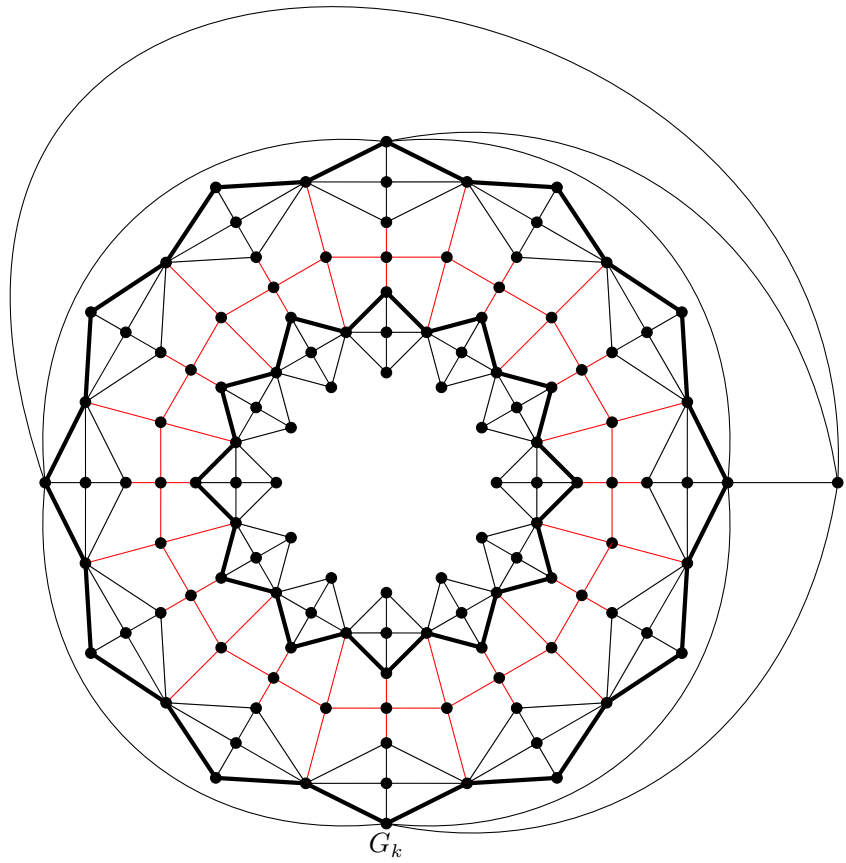

\centering
\bt
\foreach \x in {0,30,60,90,120,150,180,210,240,270,300,330}
{
\draw[color=red](\x-15:4.142)--(\x-15:3.106)--(\x-15:2.071);
\draw[color=red](\x:2.536)--(\x:3)--(\x:3.464);
\draw[color=red](\x+15:3.106)--(\x:3)--(\x-15:3.106);
\df(\x:3)circle(2pt);
\df(\x-15:3.106)circle(2pt);
\df (\x+15:4.142)--(\x:3.464)circle(2pt)--(\x-15:4.142)circle(2pt)--(\x:4)circle(2pt)--(\x:4.536)circle(2pt);
\draw[ultra thick](\x+15:4.142)--(\x:4.536)--(\x-15:4.142);
\draw (\x+15:4.142)--(\x:4)--(\x:3.464);
\df (\x+15:2.071)--(\x:1.464)circle(2pt)--(\x-15:2.071)circle(2pt)--(\x:2)circle(2pt)--(\x:2.536)circle(2pt);
\draw[ultra thick](\x+15:2.071)--(\x:2.536)--(\x-15:2.071);
\draw (\x+15:2.071)--(\x:2)--(\x:1.464);
}
\foreach \x in {0,90,180,270}
{
\draw(\x:4.536)..controls(\x+30:5.6)and(\x+60:5.6)..(\x+90:4.536);
}
\df(6,0)circle(2pt)--(4.536,0);
\draw (6,0)..controls(30:6.4)and (60:6)..(90:4.536);
\draw (6,0)..controls(-30:6.4)and (-60:6)..(-90:4.536);
\draw (6,0)..controls(50:10)and (130:12)..(-180:4.536);
\draw (0,-4.536)node[anchor=90]{$G_k$};
\et
\caption{\label{gk}$G_k$ from $H_k$}
\end{figure}
\begin{proof}
This construction comes from the motivation that $B_{5,b}$(Figure \ref{b5}) is the optimal triangular block in terms of edge-face ratio, which will be shown in later proof as we achieve $7f(B_{5,b})-4e(B_{5,b})=0$. We will do induction on the construction. The base case is $H_0$ in Figure \ref{h0}, where there are $49$ vertices and $104$ edges. The inductive step is shown in Figure \ref{hk}, where we add a layer of $288$ vertices and $672$ edges. Finally we get $G_k$ from $H_k$ by adding $1$ vertex and $8$ edges, as shown in Figure \ref{gk}. The nice thing about $G_k$ is that it only contains $B_{5,b}$ as triangular block and each boundary edge of such triangular block is adjacent to a face of length $7$. Indeed, $G_k$ has $288k+50$ vertices and $672k+112$ edges for each $k\geq 0$ and $672k+112={7\over 3}(288k+50-2)$. Clearly $G_k$ contains no $C_6$ or $K_4$ as subgraph so the proof is complete.
\end{proof}

\section{Proof of Theorem \ref{t3}}
\begin{proof}
Let $G$ be a $\{K_4,C_6\}$-free plane graph on $n$ $(n\geq 6)$ vertices with $\delta(G)\geq 3$. Denote the number of faces in $G$ as $f(G)$. By Euler's formula, $e(G)\leq {7\over 3}(n-2)$ is equivalent with $7f(G)-4e(G)\leq 0$. Observe that $7f(G)-4e(G)=\sum_{B\in B(G)}(7f(B)-4e(B))$. Hence it suffices to prove $7f(B)-4e(B)\leq 0$ for each $B\in B(G)$. We will do case analysis on different $B$'s. Since $G$ is $\{K_4,C_6\}$-free, there are only $5$ possible different types of $B$.

\begin{lemma}\label{l1}
Let $F$ be a face of length $4$ in a $\{K_4,C_6\}$-free plane graph. Then, if $F$ is adjacent to more than one non-trivial triangular blocks($B\neq K_2$), then $F$ is adjacent to exactly two non-trivial triangular blocks each being a $K_3$ that shares a vertex. These two blocks together with $F$ would look like Figure \ref{f1}:

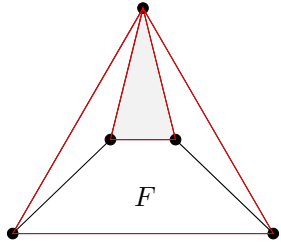
\begin{figure}
\centering
\begin{tikzpicture}
\dff (90:2)--(30:0.5)--(150:0.5)--cycle;
\df (150:0.5)circle(2pt)--(90:2)circle(2pt)--(-150:2)circle(2pt)--(-30:2)circle(2pt)--(30:0.5)circle(2pt);
\draw (-30:2)--(90:2)--(30:0.5)--(150:0.5)--(-150:2);
\draw[red] (90:2)--(30:0.5)--(150:0.5)--cycle;
\draw[red] (90:2)--(-150:2)--(-30:2)--cycle;
\draw (-0.25,-0.5)node[anchor=180]{$F$};

\end{tikzpicture}
\caption{\label{f1}When $F$ Adjacent to Two Triangular Blocks}
\end{figure}
\end{lemma}
The gray area is not a face and there are other vertices inside. Observe that all $6$ red edges have to be incident to some face of length at least $5$, thus we can uniquely concatenate the two $K_3$ triangular blocks together with $F$ to form a new ``block" with $3$ faces and $8$ edges. Let $B$ be this block, then
\begin{align*}
7f(B)-4e(B)\leq 7\cdot (3+{6\over 5})-4\cdot 8=-2.6<0.
\end{align*}
We will search and calculate the face and edge contribution of this block before we search for other triangular blocks. Later on, when we find a non-trivial triangular block adjacent to some face of length $4$, we can uniquely assign that face to the triangular block to calculate face and edge contribution.

\textbf{Case 1: }$B=K_2$.

Since $\delta(G)\geq 3$, this edge must be incident to two faces each with length at least $4$. Thus $$7f(B)-4e(B)\leq 7\cdot ({1\over 4}+{1\over 4})-4=-0.5<0.$$

\textbf{Case 2: }$B=K_3$. 

If at least one of edges of $B$ is not adjacent to some face of length $4$, we have $$7f(B)-4e(B)\leq 7\cdot (1+{1\over 5}+{2\over 4})-4\cdot 3=-0.1<0.$$On the other hand, if all three edges of $B$ is adjacent to some face of length $4$, there are  two possible drawings, as shown in Figure \ref{allf4}. Similar to other figures, the gray area is not a face and contain other vertices inside.
\begin{figure}
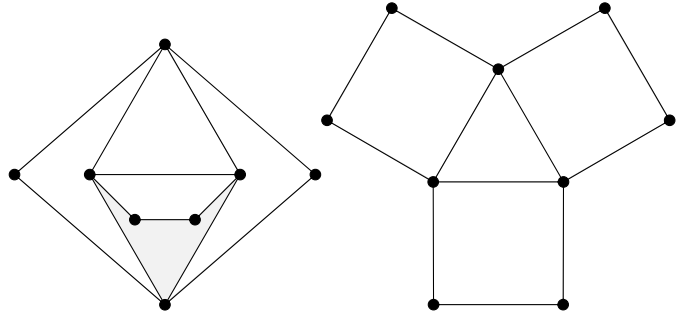

\centering
\bt
\dff(-1,0)--(0,-1.732)--(1,0)--(0.4,-0.6)--(-0.4,-0.6)--cycle;
\df(-1,0)circle(2pt)--(1,0)circle(2pt)--(0,1.732)circle(2pt)--(-2,0)circle(2pt)--(0,-1.732)circle(2pt)--(2,0)circle(2pt);
\df(-1,0)--(-0.4,-0.6)circle(2pt)--(0.4,-0.6)circle(2pt)--(1,0);
\draw(0,1.732)--(-1,0)--(0,-1.732)--(1,0);
\draw(2,0)--(0,1.732);
\et
\bt
\df(90:1)circle(2pt)--(210:1)circle(2pt)--(-30:1)circle(2pt)--(-68:2.3)circle(2pt)--(248:2.3)circle(2pt);
\draw(248:2.3)--(210:1);
\df(210:1)--(172:2.3)circle(2pt)--(128:2.3)circle(2pt)--(90:1);
\df(90:1)--(52:2.3)circle(2pt)--(8:2.3)circle(2pt)--(-30:1);
\draw(90:1)--(-30:1);
\et
\caption{\label{allf4}When all edges of $K_3$ adjacent to some face of length $4$}
\end{figure}
By Lemma \ref{l1}, we can uniquely assign the three faces of length $4$ to $B$. Then, $$7f(B)-4e(B)\leq 7\cdot (4+{9\over 4})-4\cdot 12=-4.25<0.$$

\textbf{Case 3: }$B=B_4$.

If no edges of $B$ is adjacent to some face of length $4$, then $$7f(B)-4e(B)\leq 7\cdot (2+{4\over 5})-4\cdot 5=-0.4<0.$$

If $B$ is adjacent to exactly one face of length $4$, then the only possible drawing is shown in Figure \ref{b41f4}.
\begin{figure}
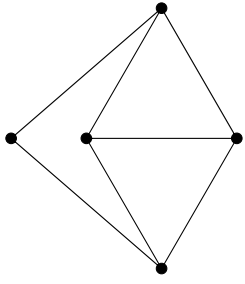

\centering
\bt

\df(-1,0)circle(2pt)--(1,0)circle(2pt)--(0,1.732)circle(2pt)--(-2,0)circle(2pt)--(0,-1.732)circle(2pt);
\draw(0,1.732)--(-1,0)--(0,-1.732)--(1,0);
\et

\caption{\label{b41f4}When $B_4$ adjacent to exactly $1$ face of length $4$}
\end{figure}
By Lemma \ref{l1}, we assign the face of length $4$ to $B$ and $$7f(B)-4e(B)\leq 7\cdot (3+{2\over 5}+{2\over 4})-4\cdot 7=-0.7<0.$$

If $B$ is adjacent to more than one face of length $4$, then the only possible drawing is shown in Figure \ref{b42f4}.
\begin{figure}
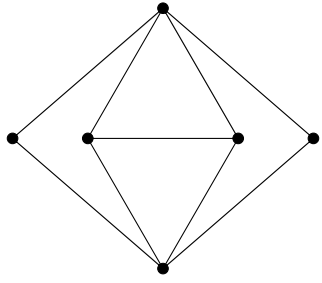

\centering
\bt

\df(-1,0)circle(2pt)--(1,0)circle(2pt)--(0,1.732)circle(2pt)--(-2,0)circle(2pt)--(0,-1.732)circle(2pt);
\df(0,1.732)--(2,0)circle(2pt)--(0,-1.732);
\draw(0,1.732)--(-1,0)--(0,-1.732)--(1,0);
\et

\caption{\label{b42f4}When $B_4$ adjacent to more than $1$ face of length $4$}
\end{figure}
By Lemma \ref{l1}, we assign the face of length $4$ to $B$ and $$7f(B)-4e(B)\leq 7\cdot (4+{4\over 4})-4\cdot 9=-1<0.$$

\textbf{Case 4: }$B=B_{5,a}$.

\begin{figure}
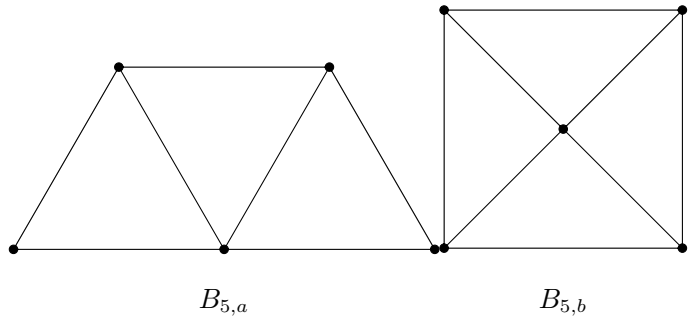

\centering
\bt[scale=0.8]
\df (0,0)circle(2pt)--(180:3.5)circle(2pt)--(120:3.5)circle(2pt)--(60:3.5)circle(2pt)--(0:3.5)circle(2pt)--(0,0);
\draw(120:3.5)--(0,0)--(60:3.5);
\draw (0,-0.5)node[below]{$B_{5,a}$};
\et
\bt[scale=0.8]
\df (45:2.8)circle(2pt)--(135:2.8)circle(2pt)--(225:2.8)circle(2pt)--(-45:2.8)circle(2pt)--(45:2.8);
\df (45:2.8)--(0,0)circle(2pt)--(135:2.8);
\draw (225:2.8)--(0,0)--(-45:2.8);
\draw (0,-2.5) node[below]{$B_{5,b}$};

\et

\caption{\label{b5}Triangular blocks on five vertices}
\end{figure}
\begin{figure}
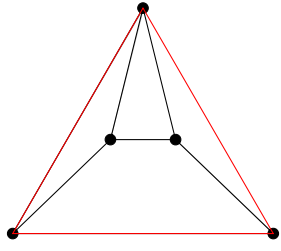

\centering
\bt

\df (150:0.5)circle(2pt)--(90:2)circle(2pt)--(-150:2)circle(2pt)--(-30:2)circle(2pt)--(30:0.5)circle(2pt);
\draw (90:2)--(30:0.5)--(150:0.5)--(-150:2);
\draw[red](-30:2)--(90:2)--(210:2)--cycle;
\et
\caption{\label{b5af4}When $B_{5,a}$ adjacent to a face of length $4$}
\end{figure}
If no edges of $B$ is adjacent to some face of length $4$, then $$7f(B)-4e(B)\leq 7\cdot (3+{5\over 5})-4\cdot 7=0.$$

Remember that we do not allow the appearance of $C_6,K_4$, or vertex of degree $2$. Thus if $B$ is adjacent to some face of length $4$, there is only one possible drawing, as shown in Figure \ref{b5af4}.

Observe that the red edges have to be adjacent to a face of length at least $7$. By Lemma \ref{l1}, we can uniquely assign the face of length $4$ to $B$ and $$7f(B)-4e(B)\leq 7\cdot (4+{3\over 7})-4\cdot 8=-1<0.$$

\textbf{Case 5: }$B=B_{5,b}$.

Remember that we do not allow the appearance of $C_6,K_4$, or vertex of degree $2$. It follows that each of the red edges has to be adjacent to some face of length at least $7$. Thus, $$7f(B)-4e(B)\leq 7\cdot (4+{4\over 7})-4\cdot 8=0.$$

\end{proof}
\section{Proof of Corollary \ref{c2}}
\begin{proof}
Let $G$ be a $\{K_4,C_6\}$-free plane graph on $n$ vertices. For simplicity of notations, in this section we denote $e(G)$ as $e$. We will show that $e(G)\leq {7\over 3}(n-2)$, which is equivalent with $7n-3e\geq 14$, or $n<9$. We will repeat the following operation:

$1)$ Delete $x$ s.t. $deg(x)\leq 2$;

When we no longer can do this operation, we get an induced subgraph $G'$ with $\delta(G')\geq 3$. Denote 	$|E(G')|=e',|V(G')|=n'$.
Each time we do operation $1)$, denote the number of deleted edge as $e_d$. Observe that $7\cdot 1-3e_d\geq 1$. Thus we get $$7n-3e\geq 7n'-3e'+1\cdot(n-n').$$

In line with usual graph theoretic terminology, we call a maximal $2$-connected subgraph a \textbf{block}. Let $b$ be the total number of blocks of $G'$. Specifically, let $b_2,b_3,b_4$, and $b_5$ denote the number of blocks of size $2,3,4$, and $5$, respectively. Let $b_6$ denote the number of blocks of size at least $6$. Then we have $b=\sum_{i=2}^6b_i$.

Let us find a lower-bound of $7n-3e-7$ for blocks of size $n$. Since $n$ is fixed, we need to plug in the largest possible $e$, which ideally would be the number of edges in a triangulated planar graph.

When $n\geq 6$, by Theorem \ref{t3} we have $7n-3e-7\geq 14-7=7$.

When $n=5$, since there is no $K_4$ as subgraph of $G$, we have $7n-3e-7\geq 7 \cdot 5-3\cdot8-7=4$ as in $B_{5,b}$.

When $n=4$, since there is no $K_4$ as subgraph of $G$, we have $7n-3e-7\geq 7 \cdot 4-3\cdot5-7=6$ as in $B_{4}$.

When $n= 3$, as in $K_{3}$ we have $7n-3e-7\geq 7\cdot3-3\cdot3-7=5$.

When $n= 2$, as in $K_{2}$ we have $7n-3e-7\geq 7\cdot2-3\cdot1-7=4$.
\\

Denote the vertex number of the $i$th block as $n_i$ and edge number as $e_i$. Using the lower bounds presented above, when $n>n'$, we have
\begin{align*}
7n'-3e'&\geq7(\sum_{i=1}^{b}n_i-(b-1))-3\sum_{i=1}^{b}e_i\\
&=\sum_{i=1}^b(7n_i-3e_i-7)+7\\
&\geq 7b_6+4b_5+6b_4+5b_3+4b_2+7\\
7n-3e&\geq 7b_6+4b_5+6b_4+5b_3+4b_2+7+(n-n').
\end{align*}

If $b_6\geq 1$, then righthand side is at least $14$, as desired.

So, let us assume that $b_6=0$. It follows that
\begin{align*}
b&=\sum_{i=2}^5b_i\\
n&\leq \sum_{i=2}^6ib_i+(n-n')\\
7n-3e&\geq 7b_6+4b_5+6b_4+5b_3+4b_2+7+(n-n')\\
&\geq 7+\sum_{i=2}^5{4\over 5}ib_i+{4\over 5}(n-n')\\
&\geq 7+{4\over 5}n.
\end{align*}
This is greater than or equal to $14$ if $n\geq 9$.

Finally, if $n=n'$ and $n\geq 2$, since the last two vertices we delete has at most one edge between them, we know $$7n-3e\geq 7\cdot 2-3\cdot 1+(n-2)=11+(n-2)$$which is greater than or or equal to $14$ when $n\geq 5$. Since $5<9$, we complete the proof of Corollary \ref{c2}.
\end{proof}

\bibliographystyle{abbrv}
\bibliography{k4c}

\end{document}